\documentclass[10pt,twoside]{article} 
\usepackage{amsmath,amsthm} 
\usepackage{amsfonts} 
\usepackage{newecpart} 
\title{Last Passage Percolation in Macroscopically Inhomogeneous Media} 
\shorttitle{LPP in Macroscopically Inhomogeneous Media} 
\numauthor{2} 
\authorone{Leonardo T. Rolla\footnote{Partially supported by CNPq and FAPERJ}} 
\authortwo{Augusto Q. Teixeira} 
 
\addressone{%
Instituto de Matem\'atica Pura e Aplicada 
Estr Dona Castorina, 110. Rio de Janeiro 22.460-320. Brazil 
} 
\addresstwo{
Instituto de Matem\'atica Pura e Aplicada 
Estr Dona Castorina, 110. Rio de Janeiro 22.460-320. Brazil 
} 
\emailone{%
leorolla@impa.br 
} 
\emailtwo{%
augusto@impa.br 
} 
 
\newcommand{\Submitted}{July 17, 2007} 
\newcommand{\Accepted}{January 14, 2008} 
\newcommand{\SubjectClassification}{%
60K35, 82D60, 60F10. 
} 
\newcommand{\Keywords}{%
last passage percolation, inhomogeneous media, variational problem, scaling limit, tasep. 
} 
\newcommand{\Abstract}{%
In this note we investigate the last passage percolation model in the presence of macroscopic inhomogeneity. 
We analyze how this affects the scaling limit of the passage time, leading to a variational problem that provides an ODE for the deterministic limiting shape of the maximal path. We obtain a sufficient analytical condition for uniqueness of the solution for the variational problem. Consequences for the totally asymmetric simple exclusion process are discussed. 
}

\usepackage{psfrag} 
\usepackage{amsmath} 
\usepackage{amssymb} 
\usepackage{amsthm} 
\usepackage{amscd} 
\usepackage{pifont}

\newcommand{\N}{\mathbb{N}} 
\newcommand{\Z}{\mathbb{Z}} 
\newcommand{\R}{\mathbb{R}}

\newtheorem{theo}{Theorem}

\begin{document} 
 
\maketitle

\section{The model and results}

The last passage percolation process has been widely studied over the last few years~%
\cite{baik01,baik05,bodineau05,hambly07,johansson00,martin04,martin06}.
There are several equivalent physical interpretations for the model.
Examples include zero-temperature directed polymer in a random environment, a certain growth process, queuing theory, a randomly increasing Young diagram and random partitions~%
\cite{glynn91,oconnell03,widom02}.
By a simple coupling argument, results obtained for last passage percolation have their duals for the totally asymmetric simple exclusion process and thus the former model is often useful for the study of the later one~%
\cite{ferrari05,ferrari06,mountford05,prahofer02,seppalainen98,seppalainen99}.
It is for the two-dimensional case that the most explicit results and estimates are known. In particular for geometric or exponential distributions, to which an exact solution was given by~\cite{johansson00}.

In this note we study last passage percolation with exponentially distributed passage times in the presence of macroscopic inhomogeneity.
We begin by restating known results with a different point of view: instead of taking the limit of a large rectangle on the usual lattice we consider the limit of a fine lattice on some fixed rectangle, which is equivalent but resembles hydrodynamics.
A continuous function $\alpha$ is defined on the macroscopic rectangle and locally modifies the parameter of the process.
The problem of studying the random microscopic path of maximal passage time leads to the variational problem of finding a deterministic macroscopic curve maximizing a certain functional.
We shall see that the rescaled passage time indeed converges to that given by the variational problem, and give sufficient analytical conditions for convergence of the maximal path's shape to a deterministic curve.

The viewpoint adopted here has immediate implications for the totally asymmetric simple exclusion process.
On the scaling limit we can describe the behavior of the total current through the origin up to a given time when the jump rate has macroscopic fluctuations in space (as well as in blocks of particles, or even both).
In particular for spatial inhomogeneity it is easy to see that the instantaneous current is non increasing in time.
The analysis of the time taken for a given amount of particles to cross the origin and of the corresponding path that gives that passage time helps understanding the bottleneck, that is, which chain of events was responsible for that delay.
The relation with TASEP will be discussed further in Section~\ref{sec:tasep}.

The last passage percolation problem may be formulated as follows. Given the origin and a point $(l,b)\in\R^2$, $|b|<l$, consider the rectangle $Q=\bigl\{(x,y):0\leqslant x\leqslant l, |y|\leqslant x, |b-y| \leqslant l-x\bigr\}$, that is, the rectangle like displayed in Figure~\ref{fig:strips} having $(0,0)$ and $(l,b)$ as vertices.
For $N\in\N$, take the grid $S_N=Q\cap\frac1N\widetilde\Z^2$, where $\widetilde\Z^2=\bigl\{(x,y)\in\Z^2:x+y\in2\Z\bigr\}$.
Let $\xi_N$ denote a field $(\xi_p)_{p\in S_N}$, whose coordinates are i.i.d. distributed as $\frac\alpha N\exp(1)$, $\alpha>0$.
One should think of $\xi_p$ as a \emph{reward} standing at the site $p$.
Let $\Pi_N$ denote the set of oriented paths $(p_0,\dots,p_{k})$ crossing $S_N$, i.e., $p_0=0$, $p_i\in S_N$, $p_i-p_{i-1}=\frac1N(1,\pm1)$ for $i=1,\dots,k$ and $p_k$ is the rightmost point of $S_N$.
If $\pi\in\Pi_N$, we define $\xi_N\cdot \pi$ as the sum of the random variables $\xi_p$ for $p\in\pi$, i.e., 
the total value of the rewards along that path.
Let $\pi^N$ denote the (a.s. unique) path that maximizes $\xi_N\cdot \pi$ and let $G(\xi_N)=\max_{\pi\in\Pi_N} (\xi_N \cdot \pi)=\xi_N \cdot \pi^N$.
It is known~\cite{johansson00} that the maximal value $G(\xi_N)$ approaches $\alpha\, l\, \gamma\bigl(\frac bl\bigr)$, where $\gamma(w) = 1 + \sqrt{1-w^2}$, as $N$ increases. Moreover, the probability of deviating $\epsilon$ from this limit decays exponentially fast in $N$.
A consequence of this fact is that the maximal path $\pi^N$ approaches the straight line connecting the origin and $(l,b)$ in the $\|\cdot\|_{\sup}$.

We now describe how the macroscopic inhomogeneity is introduced. Let $\alpha$ be a nonnegative, continuous function defined on $Q$ and for each $N\in\N$ we take $(\xi_p)_{p\in S_N}$ as independent random variables, each one distributed as $\frac{\alpha(p)}N\exp(1)$. We are interested in understanding the value of $G(\xi_N)$ and the shape of $\pi^N$ for large values of $N$.

Let us give some simple heuristic arguments.
Take some $y(\cdot)\in X$, 
where
\begin{equation*}
\label{eq:X}
 X = \bigl\{y:[0,l]\to\R \ \big|\  y(0)=0 ,y(l)=b, y \mbox{ is Lipschitz with constant } 1\bigr\}.
\end{equation*}
We first look at the macroscopic limit of $\max_\pi(\xi_N \cdot \pi)$ restricted to paths $\pi$ that stay ``pinned'' to the curve $y$.
Take $\Delta x \ll l$, consider the set $\Pi_N'\subset\Pi_N$ of paths that pass by $\bigl(x,y(x)\bigr)$ for all $x=n\Delta x$, and let $G'(\xi_N)=\max_{\pi\in\Pi_N'}(\xi_N \cdot \pi)$.
Since $\alpha$ is nearly constant on a small neighborhood of $(x,y)$,
we expect the contribution to $G'$ obtained 
between $x$ and $x+\Delta x$ to be given by $\alpha(x,y)\gamma\bigl(\frac {\Delta y}{\Delta x}\bigr)\Delta x$,
in accordance with the constant case discussed above.
This motivates the definition of the functional
\begin{equation*}
 \label{eq:G}
  \mathcal{G}(y) = \int_0^l \alpha\bigl(x,y(x)\bigr) \gamma\left(\dot y\right) dx, \quad y\in X,
\end{equation*}
where $\dot y=dy/dx$.
It is reasonable to expect that $G'(\xi_N)$ will approach $\mathcal{G}(y)$ when $N\Delta x\to\infty$ sufficiently fast and $\Delta x\to0$.
Informally we say that that the maximal path $\pi$, restricted to be pinned to the curve $y$, will asymptotically catch a total of $\mathcal{G}(y)$ in rewards. Now drop this restriction, i.e., take the maximal path among all possible paths instead of being pinned to some specific curve $y$. We guess $G(\xi_N)$ will be given by maximizing over $y$, that is,
\[
  G(\xi_N) \approx \sup_{y\in X} \mathcal{G}(y).
\]

Our first result confirms the above argument.
\begin{theo}
\label{theo:sup}
Let $\alpha$ be a continuous function on the rectangular domain $Q$. Consider the inhomogeneous last passage percolation problem as described above and take $$\mathcal{G}^* = \sup_{y\in X} \mathcal{G}(y).$$
Then $G(\xi_N) \to \mathcal{G}^*$ a.s. 
Moreover, for any $\delta>0$ there are $c,C>0$ such that
$$ P\bigl(|G(\xi_N) - \mathcal{G}^*| > \delta\bigr) < Ce^{-cN}$$
holds for all $N\in\mathbb{N}$.
\end{theo}

It is clear that the variational problem defining $\mathcal{G}^*$ is crucial for the understanding of the process. Uniqueness of its maximizer is required in order to establish convergence of the maximal paths.
Although one can easily build a function $\alpha$ for which two distinct curves maximize $\mathcal{G}^*$, the authors believe that the set of $\alpha$'s for which $y^*$ is unique is generic in the $C^0(Q)$ topology.
The next result concerns the shape of the maximal path.
\begin{theo}
\label{theo:path}
For $\alpha$ continuous there is $y^*\in X$ such that $\mathcal{G}^* = \mathcal{G}(y^*)$.
If such $y^*$ is unique, the random maximal path $\pi^N$ approaches the curve $y^*$ for $N$ large.
That is, $\| \pi^N - y_* \|_{\sup} \to 0$ a.s., and, for any $\delta>0$ there are $c,C>0$ such that
$$ P\bigl(\| \pi^N - y_* \|_{\sup} > \delta\bigr) < Ce^{-cN}$$ holds for all $N\in\mathbb{N}$.

Suppose also that $\alpha$ has continuous
derivatives $\alpha_x$ and $\alpha_y$. Then the Euler-Lagrange equation associated with the variational problem of $\mathcal{G}^*$
\begin{equation}
\label{eq:edo}
\left\{
\begin{array}{l}
\ddot y = - \frac {1} {\alpha} \bigl[ \alpha_y (1- \dot y ^2)^{3/2} + (\alpha_x \dot y + \alpha_y)(1-\dot y^2) \bigr]
\\
y(0)=0, \quad y(l)=b
\end{array}
\right.
\end{equation}
has at least one solution $y_0$.
If in addition $\alpha$ is such that $\mathcal{G}$ is strictly concave on $X$, the solution $y_0$ of~(\ref{eq:edo}) is unique and it is also the unique maximizer of $\mathcal{G}$.
\end{theo}

When $\alpha$ is constant it is easy to see that the maximizer must be a straight line, some other nontrivial examples of uniqueness are given by the condition we discuss below.
Of course uniqueness holds whenever $\mathcal{G}$ is strictly concave.
When $\alpha$ is smooth, a sufficient condition for strict concavity of $\mathcal{G}$ is
\begin{equation}
\label{eq:strongconcavity}
\alpha_{yy} < 0 \quad \mbox{and} \quad -\alpha\alpha_{yy} \geqslant \frac {1} {2} (\alpha_y)^2
\end{equation}
for all $(x,y)\in Q$.

In Section~\ref{sec:convergence} we prove Theorem~\ref{theo:sup} with the aid of large deviation estimates already known to hold when $\alpha$ is constant.
Section~\ref{sec:concave} contains the proofs of Theorem~\ref{theo:path} and the sufficiency of the condition~(\ref{eq:strongconcavity}).
In Section~\ref{sec:tasep} we discuss some implications of these results for the TASEP.

The authors believe that it is possible to prove similar results for the case of geometric distribution of passage times and for the case of non-homogeneous poissonian points with the same arguments presented in this note.

\paragraph{Acknowledgments}

We are grateful to V.~Sidoravicius for several suggestions and encouragement. We would also like to thank the anonymous referees, whose suggestions have considerably simplified and improved our exposition.

\section{Convergence of the passage time}
\label{sec:convergence}

In this section we sketch the proof of Theorem~\ref{theo:sup}. The proof consists on approximating $\alpha$ by a constant on smaller regions of the domain $Q$ and then applying the large deviation principle known to hold on these regions.

Since $\alpha$ will be approximated by other functions, we write $\mathcal{G}_\alpha(y)$ and $\xi_\alpha \cdot \pi$ to specify which $\alpha$ is being considered. Besides, processes with different $\alpha$'s can be constructed on the same probability space in a way that  $\xi_\alpha\geqslant\xi_{\tilde\alpha}$ whenever $\alpha\geqslant\tilde\alpha$.

\begin{figure}[htb]
 \hspace*{13mm}
 \includegraphics[height=47mm]{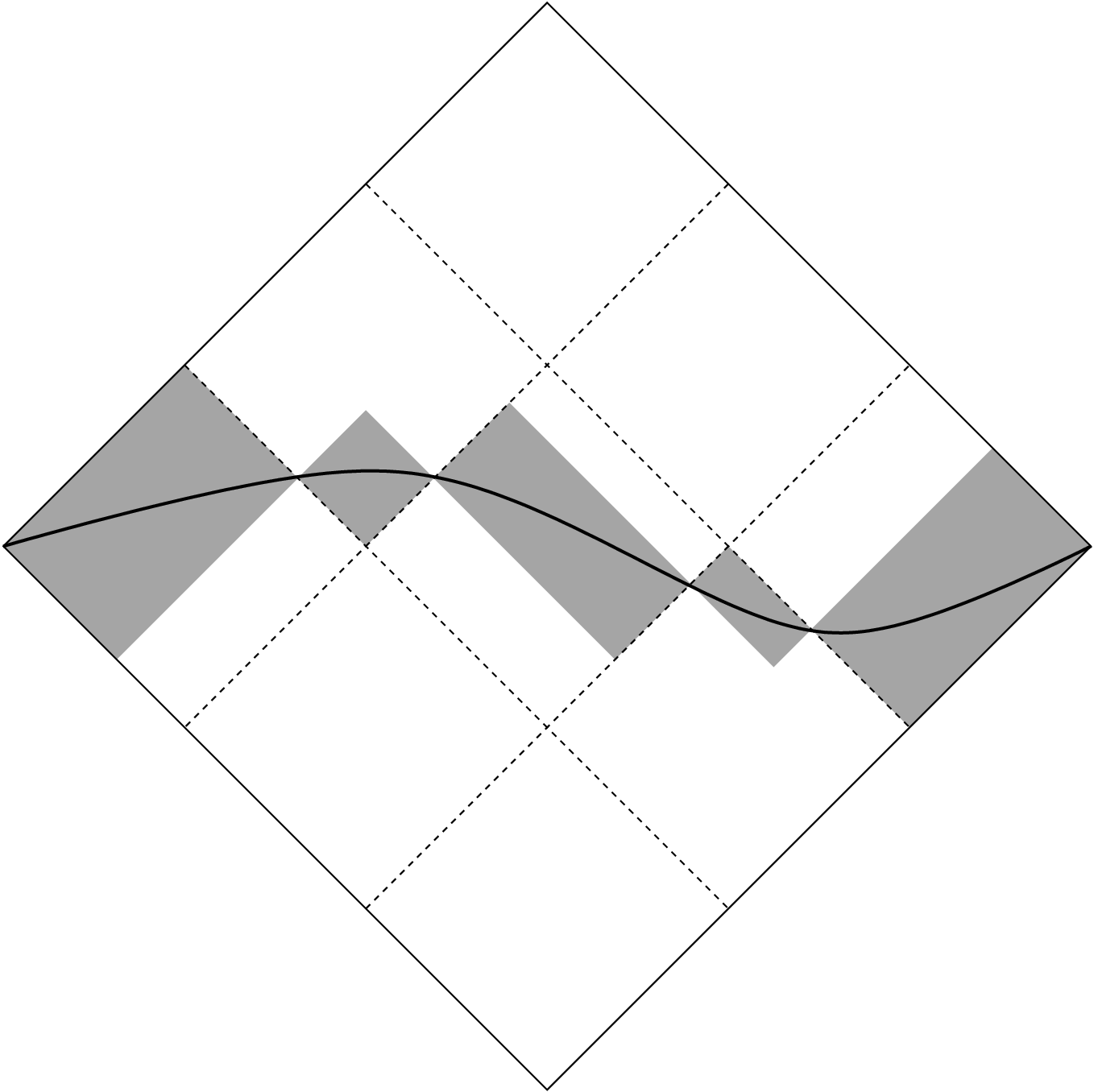}
 \hspace{20mm}
 \includegraphics[height=47mm]{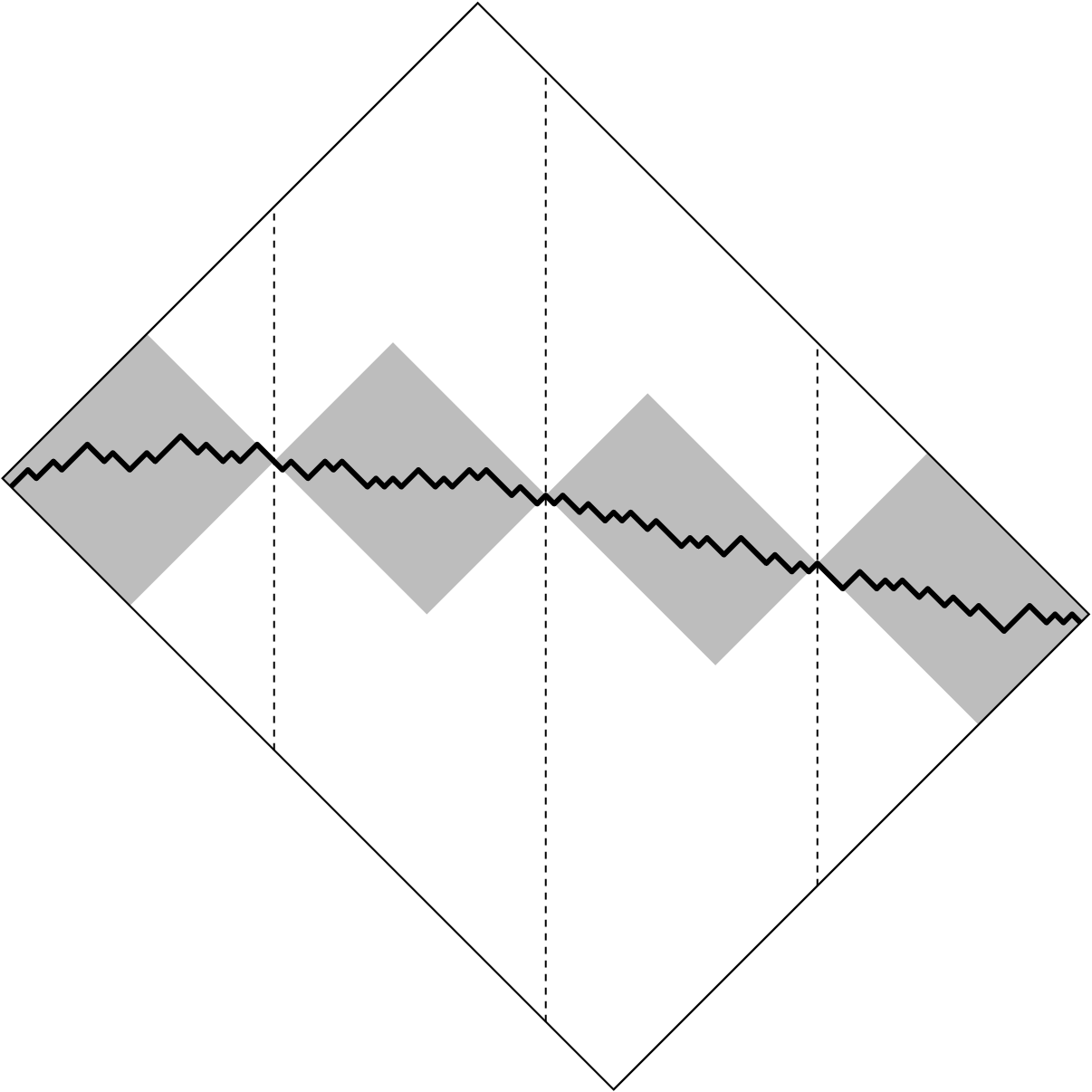}
 \caption{On the left, the definition of the sets $Q_j$ of a given $y$. 
On the right, the $k$ strips of $Q$ and a given random path $\pi^N$.}
\label{fig:Qj}
\label{fig:strips}
\end{figure}

We first prove that $P\bigl( G(\xi_N) < \mathcal{G}^* - \epsilon\bigr) < Ce^{-cN}$.
Given $\epsilon>0$, divide $Q$ as in Figure~\ref{fig:Qj}, the resulting squares being fine enough so that the value of $\alpha$ does not change more than $\frac\epsilon{6l}$ inside each of them.
Let $\tilde\alpha$ be constant on each square and given by the infimum of $\alpha$ on that square. In this case $\alpha - \frac\epsilon{6l} \leqslant \tilde{\alpha} \leqslant \alpha$.
Take $y\in X$ satisfying $\mathcal{G}(y)>\mathcal{G}^*-\epsilon/3$.
Since $\gamma\leqslant2$ we have
\[
\mathcal{G}_{\tilde\alpha}(y) \geqslant \mathcal{G}_\alpha(y) - \frac\epsilon3 \geqslant \mathcal{G}^* - \frac {2\epsilon} {3}.
\]

Let $Q_1, \dots, Q_k$ be the rectangles that contain the curve $y$ between its intersections with the grid as in Figure~\ref{fig:Qj}. $\mathcal{G}_\alpha^*(Q_j)$ stands for the supremum of $\mathcal{G}$ restricted to $Q_j$ (the functions and the integration are both restricted to this set, $Q_j$ is a rectangle and the functions connect both extremes).

Take $\Pi_N'\subseteq\Pi_N$ as the set of paths that pass through the left and the rightmost points of every $Q_j$. Given $\pi\in\Pi_N'$ write $\pi_j$ for the restriction $\pi$ to  $Q_j$.	
If $(\xi \cdot \pi^N) < \mathcal{G}^* - \epsilon$, the same must hold for every $\pi$ in $\Pi_N'$ (recall that $\pi^N$ is the maximizer of $\xi\cdot\pi$).
This implies that for some $j$, $\xi_{\tilde\alpha} \cdot \pi_j\leqslant \xi_\alpha \cdot \pi_j < \mathcal{G}^*_{\tilde \alpha}(Q_j) - \epsilon/3k$ for every $\pi$ in $\Pi_N'$. But by~\cite{seppalainen98} the probability of such event decays exponentially. So,
\[
P \bigl[ \xi_\alpha \cdot \pi^N < \mathcal{G}^* - \epsilon \bigr] \leqslant k \max_j C_j e^{-c_j N}\leqslant \tilde{C} e^{-\tilde{c} N},
\]
for suitable $\tilde{C}$ and $\tilde{c}$.
We point out that in fact $P \bigl[ \xi_\alpha \cdot \pi^N < \mathcal{G}^* - \epsilon \bigr] \leqslant Ce^{-cN^2}$, since the large deviation principle proved by~\cite{seppalainen98} is on $N^2$ and not $N$.

Finally let us show that $P\bigl( G(\xi_N) > \mathcal{G}^* + \epsilon \bigr) < Ce^{-cN}$.
Given $\epsilon>0$, take $\delta>0$ such that $|\alpha(p_1) - \alpha(p_2)|<\epsilon/8l$ whenever $|p_1-p_2|<2 \delta$. Divide $Q$ in $k$ strips of width $r = l/k$ as in Figure~\ref{fig:strips}, with $k$ such that $r < \delta$, and for each path $\pi\in\Pi_N$ write $\pi_j$ for the restriction of $\pi$ to the $j$-th strip (note that $\xi \cdot \pi=\sum_j\xi \cdot \pi_j$).

For  a given $\pi\in\Pi$, take $p_1 = (x_1,y_1),\dots,p_{k-1} = (x_{k-1},y_{k-1})$ as the points where the trajectory of $\pi$ intersects the boundary of each strip and let $y_k=b$. Note that $p_j$ needs not to be on the lattice.

We will use the following consequence of Theorem~1.6 in~\cite{johansson00}. If $\alpha \leqslant M$ is constant, there are $C,c>0$, depending only on $M$ and $\epsilon$, and such that, for any choice of $j$, $y_j$ and $y_{j+1}$,
\[
 P\left[\exists\pi_j, \xi \cdot \pi_j > \alpha r \gamma \bigl( {(y_{j+1}-y_j)}/{r} \bigr)  + \epsilon\right]
<  Ce^{-cN}.
\]

Denote by $\mathcal{G}^*_\alpha[p_1, p_2]$ the supremum of $\mathcal{G}_\alpha$ taken over the $1$-Lipschitz functions that connect $p_1$ with $p_2$, the integration done on the appropriate domain.
Now
\[
\mathcal{G}^*_\alpha \geqslant \sum_{j=0}^{k-1} \mathcal{G}^*_\alpha[p_j,p_{j+1}] \geqslant
\sum_{j=0}^{k-1} \left[ \alpha(p_j) + \frac\epsilon{8l} \right] r \gamma \bigl((y_{j+1}-y_j)/r \bigr) - \frac\epsilon2.
\]
The first inequality is obvious. The second one follows by the choice of $\delta$, since it implies that for any 1-Lipschitz function $y$ connecting $p_j$ and $p_{j+1}$ we have $\alpha\bigl(x,y(x)\bigr) \geqslant \alpha(p_j) - \frac \epsilon {8l}$ and $\gamma\leqslant2$.

Therefore, if $[\xi_\alpha \cdot \pi^N]>\mathcal{G}^*_\alpha + \epsilon$, there must be $j$ such that $\xi_{\bar\alpha} \cdot \pi_j\geqslant\xi_\alpha \cdot \pi_j \geqslant 
\bar\alpha r \gamma \bigl((y_{j+1}-y_j)/r \bigr) + \frac\epsilon{2k}$,
where $\bar\alpha \equiv \alpha(p_j) + \frac\epsilon{8l}$ satisfies $\bar\alpha \geq \alpha(p)$ for all $p\in\pi_j$. Considering all possible choices of $j$, $y_j$, $y_{j+1}$ we get
\[
\textstyle P\bigl[G_\alpha(\xi) > \mathcal{G}^*_\alpha + \epsilon \bigr] \leqslant 16rlkN^2Ce^{-cN} \leqslant \tilde Ce^{- \tilde cN},
\]
concluding the proof of Theorem~\ref{theo:sup}.

\section{Limiting shape of the maximal path}
\label{sec:concave}

In this section we present the proofs of Theorem~\ref{theo:path} and the sufficiency of~(\ref{eq:strongconcavity}) for uniqueness.

We first address the question of existence of a maximizer $y^*$.
By Ascoli-Arzel\'a, $X$ is a compact space with the uniform topology.
However, $\mathcal{G}$ is not a continuous functional in that topology (notice that $\mathcal{G}(y)$ is defined in terms of the derivative of $y$, which exists a.e. in $[0,l]$ for $y\in X$).
Nevertheless the proof will follow by compactness, since the $\mathcal{G}$ is upper semi-continuous, as we prove now.

We shall show that $\mathcal{G}=\inf_m \mathcal{G}_m$, where the ``Riemann sums''
\[
\mathcal{G}_m(y) = 
\sum_{i=1}^m 
\frac{l}{m} 
\left[ 
 \textstyle \sup_{\left[\frac{i-1}{m} l , \frac{i}{m} l\right]} \alpha\bigl(y(x), x\bigr)
\right]
\gamma\left({\textstyle \frac ml\bigl[y(\frac{i}{m}l) - y(\frac{i-1}{m}l)\bigr]}\right)
\]
are clearly continuous in the uniform topology if $\alpha$ is continuous on $Q$.
Since $\gamma$ is a concave function, it follows from Jensen's inequality that $\mathcal{G} \leqslant \mathcal{G}_m$ for every $m$, so it remains to show that $\mathcal{G}_m(y)\to\mathcal{G}(y)$ for all $y$.
Let $y\in X$ and take a sequence of $y_n\in C^1$ such that $||y-y_n||_\infty\to0$ and $||\dot y - \dot y_n||_{L^1}\to0$ as
$n\to\infty$.
Write
\begin{align*}
|\mathcal{G}(y)-\mathcal{G}_m (y)| & \leqslant |\mathcal{G}(y)-\mathcal{G} (y_n)| + |\mathcal{G}(y_n)-\mathcal{G}_m (y_n)| + |\mathcal{G}_m(y_n)-\mathcal{G}_m (y)|.
\end{align*}
The first term vanishes as $n\to\infty$ by the bounded convergence theorem.
The second term vanishes as $m\to\infty$ for $n$ fixed, since for smooth $y$ the numbers $\mathcal{G}_m(y)$ indeed correspond to Riemann sums of a continuous integrand along the path $x\mapsto(x,y,\dot y)$.
It thus suffices to show that the last term goes to zero uniformly in $m$ as $n\to\infty$.
Since $|\gamma(w)-\gamma(w')|\leqslant C\sqrt{|w-w'|}$ for all $w$, $w'$, 
$$ |\mathcal{G}_m(y_n)-\mathcal{G}_m (y)|
 \leqslant C \bar \alpha \sum_{i=1}^m \frac {l}{m} \bigg(\int_{i-1}^{i}|\dot y(lz/m) - \dot y_n (lz/m) |dz \bigg)^{1/2}
 \leqslant C \bar \alpha \sqrt{l ||y-y_n||_{L^1}},
$$
where the second inequality is Jensen's inequality for the square root and $\bar\alpha$ denotes $\sup_Q|\alpha|$.
This proves that $\mathcal{G}$ is upper semi-continuous,
therefore there is a $y^*$ attaining $\sup_{y\in X} \mathcal{G}(y)$.

We now discuss the convergence of $\pi^N$ to $y^*$ when this is the unique maximizer. It follows from uniqueness of $y^*$ and upper semi-continuity that
\begin{equation*}
 \mbox{for any } \delta>0,\  \exists \epsilon>0 \mbox{ such that } y\in X,\|y-y^*\|_{\sup} > \delta \mbox{ implies } \mathcal{G}(y)<\mathcal{G}^*-\epsilon.
\end{equation*}
This, together with the exponentially fast convergence in probability given by Theorem~\ref{theo:sup}, implies the desired result.

Suppose $\alpha$ has continuous derivatives and consider the Euler-Lagrange problem~(\ref{eq:edo}). For a given $C^2$ function $y\in X$ with $|\dot y|<1$, it is easy to see that~(\ref{eq:edo}) holds if and only if
$$
\frac {\partial} {\partial t} \mathcal{G}\bigl(y + t h\bigr) _{\bigl|_{t=0}} = 0
$$
for every smooth perturbation $h$ that vanishes on the extremes of the interval.
Rewrite~(\ref{eq:edo}) as
\begin{equation}
\label{eq:edo2}
\left\{
\begin{array}{l}
\dot w = - \frac {1} {\alpha} \bigl[ \alpha_y (1- w ^2)^{3/2} + (\alpha_x w + \alpha_y)(1-w^2) \bigr] \\
 \dot y = w
 \end{array} \right.
\end{equation}
and consider the initial value problem given by $y(0) = 0$, $w(0) = w_0$ and~(\ref{eq:edo2}). Existence, uniqueness and continuous dependence of solutions on the initial condition $w_0$ follow from basic ODE theory.
Now notice that for $w_0 = \pm 1$ we have $w(x)=\pm1\ \forall x$ and $y(l) = \pm l$.
So, for given $b\in(-l,l)$, by the intermediate value theorem there is $w_0\in(-1,1)$ for which the corresponding solution $y_0$ satisfies $y_0(l)=b$.
Since $\{w=\pm1\}$ is an invariant set, $y_0$ must satisfy $|\dot y_0(x)|<1\ \forall\ x$.
Therefore there is at least one solution $y_0\in X$ of~(\ref{eq:edo}).

Finally we 
show that when $\mathcal{G}$ is strictly concave, the solution $y_0$ is unique and it is also the unique maximizer of $\mathcal{G}$.
We simply adapt a standard method for showing that a critical point of concave functional must be the unique maximizer.
Let $y_0\in X$ be a solution of~(\ref{eq:edo}).
Suppose there were $y_1\in X\backslash\{y_0\}$ with $\mathcal{G}(y_1)\geqslant \mathcal{G}(y_0)$.
By the strict concavity, there is $y_2 \in X$ such that $\mathcal{G}(y_2) > \mathcal{G}(y_0)$.
Approximate $y_2$ by some smooth $y_3 \in X$ such that still $\mathcal{G}(y_3) > \mathcal{G}(y_0)$.
The concavity of $\mathcal{G}$ implies $\frac {\partial} {\partial t}  \mathcal{G}(y_0 + th)_{|_{t=0}} > 0$,
where $h=y_3-y_0$, a contradiction.
Therefore $y_0$ is the unique maximizer of $\mathcal{G}$ in $X$.
As $y_0$ was taken as any solution of the contour value problem, this must be the unique solution.

This finishes the proof of Theorem~\ref{theo:path}.
We end this section showing that~(\ref{eq:strongconcavity}) is sufficient for $\mathcal{G}$ to be strictly concave.

Let $x_0 \in (0,l)$ be fixed and define
$\alpha(y) = \alpha (x_0,y), \ z(y,w) = \alpha(y) \gamma(w)$. 
What we shall actually prove is that the function $z$ is strictly concave, which is much stronger than the functional $\mathcal{G}$ having this property.
We study the eigenvalues of the Hessian $H_z$ of $z$ to have conditions under which $z$ is a concave function.
$$
{\rm det}(H_z - \lambda) = \lambda^2 - \underset {s} {\underbrace {(\alpha \gamma'' + \alpha'' \gamma)}} \lambda + \underset {p} {\underbrace {[\alpha \alpha'' \gamma \gamma'' - (\alpha' \gamma') ^2]}}.
$$
Both eigenvalues are negative if and only if $s<0$ and $p>0$. Since
$$
\sup_w \frac {(\gamma')^2} {- \gamma \gamma''} < \frac12,
$$
it is enough to require~(\ref{eq:strongconcavity}) to have the strict concavity of $z$ and consequently of $\mathcal{G}$. 

\section{Applications to TASEP%
}
\label{sec:tasep}

The last passage percolation model can be coupled with the totally asymmetric simple exclusion process with the initial condition that the sites are occupied iff they lie to the left of the origin. See~\cite{seppalainen01} for a description of such correspondence.
The current through the origin for the TASEP can be understood by considering last passage percolation between $(0,0)$ and $(l,0)$, since this gives the rescaled time needed for $lN/2$ particles to cross the origin. At the scaling limit this converges to $\mathcal{G}^*[l]$.

The model considered in the previous sections can describe, for example, the TASEP with the jump rates depending smoothly on the position of each site, when we take $\alpha$ depending on $y$.
In the hydrodynamic limit, the instantaneous current through the origin will be a decreasing function of the time.
Intuitively, as the time passes, the process finds more bottlenecks that were not important before, and at some point such barriers start being determinant for the passage time, at least until the system meets even narrower bottlenecks.

\begin{figure}[h]
\psfrag{yaxis}{$y$}
\psfrag{xaxis}{$x$}
\psfrag{alphaaxis}{$\alpha(y)$}
\begin{center}
\includegraphics[height=4cm]{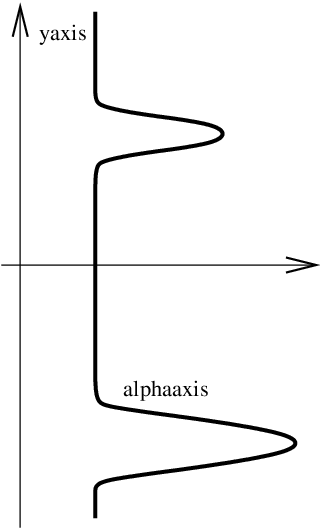}
\hspace{10mm}
\includegraphics[height=4cm]{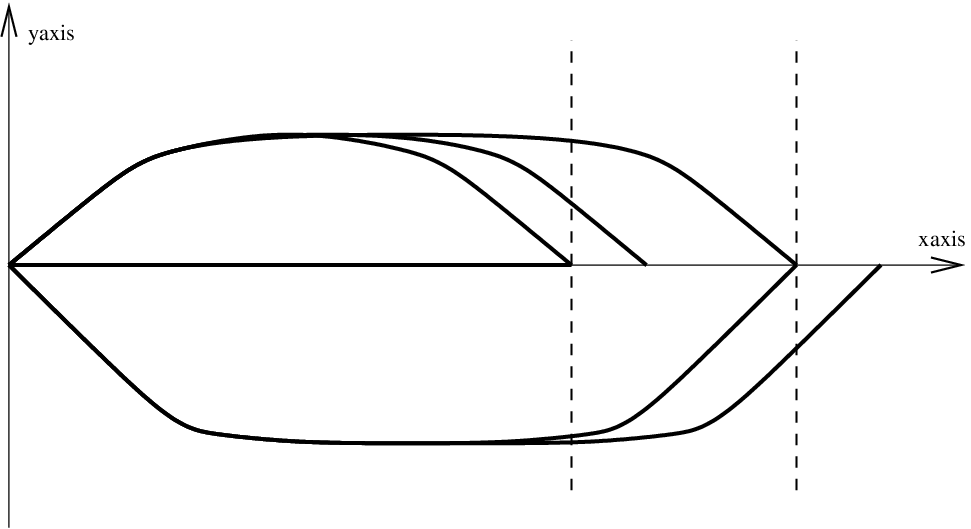}
\caption{on the left an example of $\alpha(y)$ having two strong peaks; on the right the maximizers $y^*_l$ for $l$ in three different regions.}
\label{fig:tasep1}
\end{center}
\end{figure}

The above fact follows easily from the variational formulation of $\mathcal{G}$.
%
%
Consider for instance the graph of $\alpha(y)$ having two peaks as in Figure~\ref{fig:tasep1}. For short times the system dynamics does not feel the existence of such regions of large average jump time, so $\mathcal{G}^*[l]$ (i.e., the time needed for $Nl/2$ particles to cross the origin, at the scaling limit) grows linearly with $l$ as in the homogeneous case.
For larger times the first peak will become important: there will be a traffic jam before this point and the rescaled time $\mathcal{G}^*[l]$ will increase roughly linearly with $l$, but at bigger rate, given by the value of $\alpha$ at this peak.
Now for even larger times again the system finds a harder difficulty to overcome, which is to wait for particles to cross the stronger bottleneck given by the second peak, leading to a low density of particles flowing through the origin, so from this time on $\mathcal{G}^*[l]$ increases at an even higher rate, again given by the value of $\alpha$ at the new peak.
The curves shown in Figure~\ref{fig:tasep1} give the maximizers $y^*_l$ of $\mathcal{G}^*[l]$, for $l$ in these three regions, illustrating the decreasing of instantaneous current.

For a proof of this phenomenon notice that the instantaneous current being non increasing is equivalent to saying that the amount of time needed for a certain flux to pass by the origin is a convex function. Now the following inequality implies that $l\mapsto\mathcal{G}^*[l]$ is convex. Given $l_0>0$,
$$
\mathcal{G}^*[l] \geqslant
\mathcal{G}^*[l_0]+ 2\bar\alpha (l-l_0),  
$$
where $\bar\alpha=\textstyle\sup_{[0,l_0]}\alpha\bigl(\bar y(x)\bigr)=\alpha\bigl(\bar y(x_0)\bigr)$ and $\mathcal{G}(\bar y) = \mathcal{G}^*[l_0]$.
To see why this inequality hods for $l>l_0$ we translate the curve $\bar y$ by $l-l_0$ after $x_0$ and fill the gap with a constant line; for the resulting curve $\tilde y$ we have $\mathcal{G}^*[l]\geqslant\mathcal{G}(\tilde y)=\mathcal{G}(\bar y)+ 2\bar\alpha (l-l_0)=\mathcal{G}^*[l_0]+ 2\bar\alpha (l-l_0)$.
For $l<l_0$ we find a point $x_1$ such that $\bar y(x_1) = \bar y(x_1 + l_0 - l)$ and concatenate $\bar y|_{[0,x_1]}$ with $\bar y|_{[x_1 + l_0 - l, l_0]}$.

\end{document}